\documentclass[12pt,a4paper]{article}

\usepackage[]{algorithm2e}
\usepackage[utf8]{inputenc}
\usepackage[english]{babel}
\usepackage{color}
\usepackage{amsmath}
\usepackage{amsthm}
\usepackage{amsfonts}
\usepackage{amssymb}
\usepackage{graphicx}
\usepackage{mathtools}
\usepackage[left=2cm,right=2cm,top=2cm,bottom=2cm]{geometry}
\usepackage{hyperref}
\newcommand{\Init}{\mathrm{Init}}
\newcommand{\Unsafe}{\mathrm{Unsafe}}
\newcommand{\RR}{\mathbb{R}}
\newcommand{\LL}{\mathcal{L}}
\newcommand{\XX}{\mathbb{X}}

\newtheorem*{problem}{Problem}
\newtheorem{lemma}{Lemma}
\title{
Solving  Underdetermined Boundary Value Problems By Sequential Quadratic Programming}
\author{Jan Ku\v{r}\'{a}tko\thanks{ORCID: \href{http://orcid.org/0000-0003-0105-6527}{0000-0003-0105-6527}; Faculty of Mathematics and Physics, Charles University, Czech Republic; Institute of Computer Science, The Czech Academy of Sciences} and Stefan Ratschan\thanks{ORCID: \href{http://orcid.org/0000-0003-1710-1513}{0000-0003-1710-1513}; Institute of Computer Science, The Czech Academy of Sciences}}
\begin{document}
\maketitle
\abstract{Ordinary differential equations that model technical systems often contain states, that are considered dangerous for the system. A trajectory that reaches such a state usually indicates a flaw in the design. In this paper, we present and study the properties of an algorithm for finding such trajectories. That is, for a given ordinary differential equation, the algorithm finds a trajectory that originates in one set of states and reaches another one. The algorithm is based on sequential quadratic programming applied to a regularized optimization problem obtained by multiple shooting.}
\section{Introduction}
\label{sec:Introduction}
In this paper we present and study the properties of an algorithm that, for a given ordinary differential equations, finds a trajectory of arbitrary length that originates in a set of initial states and reaches a set of unsafe states. We call such a trajectory an error trajectory. We consider both initial and unsafe sets to be ellipsoids. Unlike classical boundary value problems, this problem is underdetermined, and hence usually does not have a unique solution.

Applications in systems verification have recently motivated a lot of research on this topic, 
sometimes called the problem of ``falsification'', as the dual of ``verification''. There are several existing approaches that reduce the problem to an optimization problem~\cite{Abbas:2011,S-TaLiRo:2011,KuratkoRatschan:2014,Lamiraux:2004,Nghiem:2010,Zutshi:2013}. However, up to now, optimization has been used as a blackbox, and the algorithms have been designed without knowledge of the inner workings of optimization algorithms. In this paper, we open up the black box, exploring the specific structure of the optimization problem coming from  multiple-shooting formulations. More specifically, we compare several alternative formulations as optimization problems. We also identify the necessity of regularizing the resulting optimization problem, and study several alternative regularization terms.  The whole approach is based  on sequential quadratic programming.

Multiple shooting has been applied to this problem before, even in the more general case of hybrid dynamical systems~\cite{Zutshi:2013,KuratkoRatschan:2014}. 
However, these approaches use optimization as a black box, and use only one ad hoc formulation without any attempt for regularization of the problem. 




The contribution of this paper is the following:
\begin{itemize}
\item We formulate the underdetermined BVP problem which arises into an optimization task featuring regularization terms to mitigate the problem of having infinitely many solutions.
\item We consider various formulations of the resulting optimization problem and compare them.
\item We apply different approximation schemes for the Hessian of the Lagrangian and compare them. Especially, we are interested in sparsity preserving approximations.
\item Finally, we briefly discuss the choice of the solution technique for the saddle-point matrix in dependence on the dimension of the problem.
\end{itemize}
The structure of the paper is as follows: We formulate the problem  we try to solve and state its relation to classical boundary value problems in Section \ref{sec:ProblemFormulation}. In Section \ref{sec:Nonlinear} we reduce the problem to non-linear constrained optimization using multiple shooting. In Section~\ref{sec:OFaC} we introduce various alternative formulations as optimization problems and different regularization terms. In Section~\ref{sec:SQP}, we briefly review the sequential quadratic programming method. In Section~\ref{sec:Dis_OFaC} we study properties of the alternative optimization formulations. In Section~\ref{sec:Numerics}, we discuss practical considerations arising when implementing the resulting method. In Section~\ref{sec:Benchmarks} we present computational experiments, and in Section~\ref{sec:Conclusion} we conclude the paper.

The research published in this paper was supported by GA{\v C}R grant 15-14484S and by the long-term strategic development financing of the Institute of Computer Science (RVO:67985807).
\section{Problem Formulation}
\label{sec:ProblemFormulation}
In this section we introduce the problem we try to solve. We start with a differential equation of the form
\begin{equation}
\label{eq:DiffEq}
\frac{dx(t)}{dt} = f(t, x(t)),
\end{equation} 
where $x: \RR \to \RR^n$ is an unknown function of a variable $t \geq t_0$, and $f: \RR\times\RR^n\to\RR^n$ is  continuously differentiable. 

To stress the dependence of the solution of \eqref{eq:DiffEq} on the initial value $x_0 \in \RR^n$ we introduce the \emph{flow} function $\Phi : \RR\times \RR^n \to \RR^n$. If we fix the initial value $x_0$, then the resulting function $\Phi: \RR \to \RR^n$ expresses the solution $x(t)$ of \eqref{eq:DiffEq}. Therefore, for the initial time $t_0$ we have $\Phi(t_0, x_0) = x_0$, and for $t \geq t_0$, $\Phi(t, x_0) = x(t)$ .

 We formulate the problem we try to solve in the following way.
 \begin{problem}
 \label{prob:Fals}
 Assume a dynamical system whose dynamics is governed by the differential equations in \eqref{eq:DiffEq}. Let $\Init$ and $\Unsafe$ be sets of states in $\RR^n$. Find a trajectory of the dynamical system that starts in $\Init$ and reaches $\Unsafe$. Formally, we look for an $x_0\in\RR^n$ and $t_1\in\RR^{\geq 0}$ such that $x_0\in\Init$, and $\Phi(t_1, x_0)\in \Unsafe$.
\end{problem}

We call such a trajectory from $\Init$ to $\Unsafe$ an \emph{error} trajectory of the system. We assume that there exists an error trajectory and that the sets $\Init$ and $\Unsafe$ are disjoint. In addition, we assume the sets $\Init$ and $\Unsafe$ to be ellipsoids with centres $c_I$ and $c_U$, that is
\begin{align*}
\Init & = \left\{ v \in \RR^n \mid \left( v - c_I\right)^T E_I( v- c_I) \leq 1\right\}, \\
\Unsafe & = \left\{ v \in \RR^n \mid \left( v - c_U\right)^T E_U( v- c_U) \leq 1\right\}.
\end{align*}
We denote the norms induced by symmetric definite matrices $E_I \in \RR^{n \times n}$ and $E_U \in \RR^{n \times n}$ by $\| \cdot \|_{E_I}$, and similarly by $\| \cdot \|_{E_U}$.

Note that the problem is a BVP  with separated boundary value conditions, however, it is not in standard form~\cite[Ch.~6]{Ascher:1998}:
\begin{itemize}
\item The upper bound $t_1$ on time $t \geq t_0$ is unknown.
\item The boundary conditions are of the form
\[
g\left(x_0, \Phi(t_1, x_0)\right) =
\begin{bmatrix}
\| x_0 - c_I \|_{E_I} - 1 \\
\| \Phi(t_1, x_0) - c_U \|_{E_U} - 1
\end{bmatrix}
\leq
0,
\]
therefore, $g: \RR^{2n} \to \RR^2$.
\end{itemize}
The unknown upper bound on time $t_1$ can be eliminated by transforming the BVP into an equivalent one with a fixed upper bound, introducing one more variable~\cite[Ch.~11]{Ascher:1995}, \cite[Ch.~6]{Ascher:1998}. However, the problem we try to solve remains underdetermined.
\section{Non-linear Minimization Formulation}
\label{sec:Nonlinear}
\begin{figure}[htb]
\centering
\includegraphics[scale=1.5]{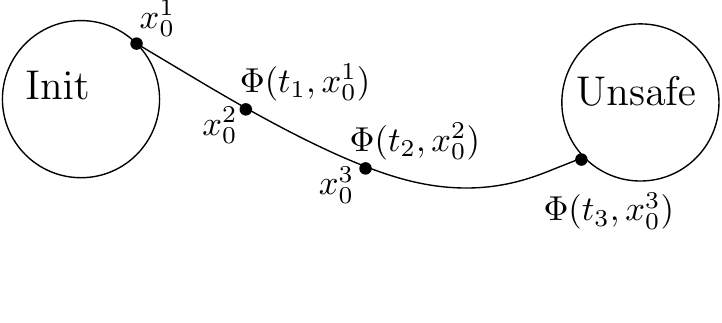}
\caption{Illustration of Multiple Shooting}
\label{fig:Traj}
\end{figure}

In this section we reformulate the problem as a minimization problem. Since we seek a trajectory which originates in one set and reaches another set we may apply techniques that are used for BVPs, in our case multiple shooting.
We solve our problem by connecting several trajectories of the system into one trajectory that is an error trajectory. Figure~\ref{fig:Traj} illustrates the overall idea on three trajectories with initial states $x_0^i$, $1 \leq i \leq 3$. The initial state $x_0^1 \in \Init$ and the final state $\Phi(t_3, x_0^3) \in \Unsafe$.

For the sequence of $N$ trajectories of lengths $t_i$,  $1 \leq i \leq N$, we define the vector of parameters
\begin{equation}
\label{eq:parameters}
\XX = \left[ x_0^1, t_1, x_0^2, t_2, \ldots, x_0^N, t_N \right]^T \in \RR^{N(n+1)}\ .
\end{equation}
We reformulate our problem so that we solve
\begin{equation}
\label{eq:minF(X)}
\min_{\XX \in \Gamma} F(\XX)\ ,
\end{equation}
where $\Gamma$ is a feasible set defined by a system of equations and inequations and $F: \RR^{N(n+1)} \rightarrow \RR$ is an objective function that measures the closeness to being an error trajectory. Note that we need regularization: In Figure \ref{fig:Traj}, if we fix $x_0^1$ and $\Phi(t_3, x_0^3)$, and constrain $x_0^2=\Phi(t_{1}, x_0^{1})$, $x_0^3=\Phi(t_{2}, x_0^{2})$, there might still be infinitely many possibilities how to choose the lengths of the segments. We consider several formulations of the objective function $F$ and the feasible set $\Gamma$ in \eqref{eq:minF(X)}.

\section{Objective Function and Constraints}
\label{sec:OFaC}
\begin{figure}
\centering
\includegraphics[scale=1.5]{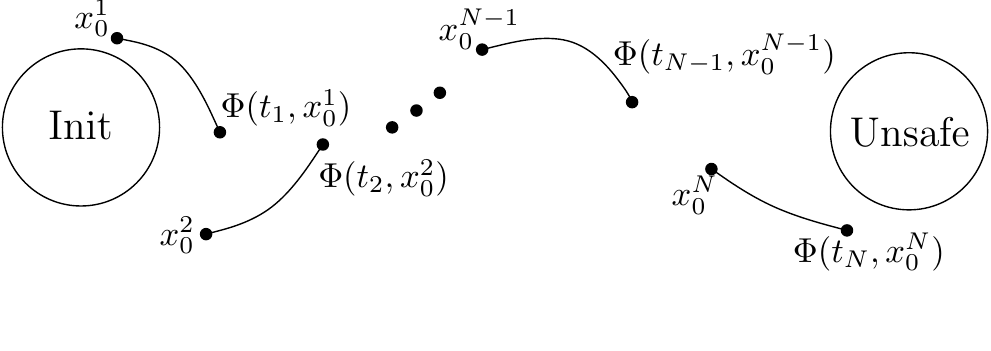}
\caption{$N$ trajectories that we aim to connect into one error trajectory.}
\label{fig:connect_seg}
\end{figure}

In order to arrive at an error trajectory from the $N$ trajectories in Figure \ref{fig:connect_seg} we need to have $x_0^1 \in \Init$ and $\Phi(t_N, x_0^N) \in \Unsafe$. In addition we need to satisfy the matching conditions $x_0^{i+1} = \Phi(t_i, x_0^i)$ for $1 \leq i \leq N-1$. For the minimization problem in \eqref{eq:minF(X)} we may formulate these  either as penalties or constraints. 
To this end we consider several different formulations of the minimization problem. We may define the objective function $F(\XX)$ in one of the following ways:
\begin{itemize}
\item $F_1(\XX) = \frac{1}{2}\left[\|x_0^1 - c_I \|_{E_I}^2 + \|\Phi(t_N, x_0^N) - c_U \|_{E_U}^2\right]$,
\begin{itemize}
\item Minimize the sum of squares of distances of $x_0^1$ from the ellipsoid centre $c_I$, and $\Phi(t_N, x_0^N)$ from the ellipsoid centre $c_U$, respectively. If an error trajectory originates at $c_I$ and reaches $c_U$, then we get the value of the objective function in the minimum equal to zero.
\end{itemize}
\item $F_2(\XX) = \frac{1}{2}\sum_{i=1}^{N-1}\|  x_0^{i+1} - \Phi(t_i, x_0^i) \|_2^2$,
\begin{itemize}
\item Minimize the sum of squares of distances between the final state and the initial state of two consecutive solutions segments. Since an error trajectory is continuous, that is $x_0^{i+1} = \Phi(t_i, x_0^i)$ for $1 < i < N$, we expect the value of the objective function to be equal to zero at the end of computation.
\end{itemize}
\item $F_3(\XX) = \frac{1}{2}\left[\|x_0^1 - c_I \|_{E_I}^2 + \sum_{i=1}^{N-1} \|  x_0^{i+1} - \Phi(t_i, x_0^i) \|_2^2 + \|\Phi(t_N, x_0^N) - c_U \|_{E_U}^2\right]$.
\begin{itemize}
\item This combines both objective functions described above. However, when an error trajectory does not originate at the centre $c_I$ and reach $c_U$, then the value of the objective function is non-zero and even the term corresponding to $F_2(\XX)$ may be non-zero.
\end{itemize}
\end{itemize}

We formulate three alternative vectors of constraints that define the feasible set $\Gamma$ in \eqref{eq:minF(X)}:
\begin{align*}
c_1(\XX) & = 
\begin{bmatrix}
x_0^2 - \Phi(t_1, x_0^1) \\
x_0^3 - \Phi(t_2, x_0^2) \\
\vdots \\
x_0^N - \Phi(t_{N-1}, x_0^{N-1})
\end{bmatrix} \ ,& 
c_2(\XX) & = 
\begin{bmatrix}
x_0^2 - \Phi(t_1, x_0^1) \\
x_0^3 - \Phi(t_2, x_0^2) \\
\vdots \\
x_0^N - \Phi(t_{N-1}, x_0^{N-1}) \\
\frac{1}{2} \left( \|x_0^1 - c_I \|_{E_I}^2 - 1\right) \\
\frac{1}{2} \left( \|\Phi(t_N, x_0^N) - c_U \|_{E_U}^2 - 1 \right) \\
\end{bmatrix} \ , \\ 
c_3(\XX) & = 
\begin{bmatrix}
\frac{1}{2} \left( \|x_0^1 - c_I \|_{E_I}^2 - 1\right) \\
\frac{1}{2} \left( \|\Phi(t_N, x_0^N) - c_U \|_{E_U}^2 - 1 \right)
\end{bmatrix}  \ ,
\end{align*}
where $c_1(\XX) \in \RR^{n(N-1)}$, $c_2(\XX) \in \RR^{n(N-1)+2}$ and $c_3(\XX) \in \RR^{2}$.

We will investigate three formulations of the optimization problem \eqref{eq:minF(X)}:
The first formulation sets penalty terms in $F(\XX)$ for distances of the initial and the final state to sets $\Init$ and $\Unsafe$, and states the matching conditions as constraints:
\begin{align}
\label{eq:pen_match}
\min_{\XX}F_1(\XX) & \quad \mathrm{subject\ to}  \quad c_1(\XX) = 0\, . 
\end{align}
The second formulation sets penalty terms for the matching conditions, and states $x_0^1 \in \Init$ with $\Phi(t_N, x_0^N) \in \Unsafe$ as constraints:
\begin{align}
\label{eq:gap_pen}
\min_{\XX}F_2(\XX) & \quad \mathrm{subject\ to}  \quad c_3(\XX) = 0 \, .
\end{align}

In addition, we will consider the unconstrained version  in which we use only penalty terms for matching conditions as well as for distances:
\begin{align}
\label{eq:unconst}
\min_{\XX}F_3(\XX)\, .
\end{align}
Finally, we may consider solving the vector of constraints $c_2(\XX) = 0$ with no objective function. This would lead to finding a solution to a system of $(N-1)n + 2$ non-linear equations with $N(n+1)$ unknowns, where $N$ is the number of segments and $n$ is the state-space dimension. When $N + n - 2 > 0$ the system is underdetermined. 

We may observe that problem formulations \eqref{eq:pen_match}--\eqref{eq:unconst} allow for infinitely many solutions. Even for a fixed length of an error trajectory there are still free parameters $t_1, \ldots, t_N$. Because of this non-uniqueness we introduce a regularization term into \eqref{eq:pen_match}--\eqref{eq:unconst} so that we can control the lengths of segments and put preference to some solutions over the others. To this end we consider the following regularization terms:
\begin{itemize}
\item $R_1(\XX) = \frac{1}{2}\sum_{i=1}^Nt_i^2$,
\begin{itemize}
\item This regularization term aims at keeping lengths $t_i$, $1 \leq i \leq N$, to be equally distributed and the length of an error trajectory to be shortest.
\end{itemize}
\item $R_2(\XX) = \frac{1}{2}\sum_{i=1}^{N-1}(t_{i+1} - t_i)^2$,
\begin{itemize}
\item By using this regularization term we try to have lengths of two consecutive segments the same, hence the length of an error trajectory gets equally distributed.
\end{itemize}
\item $R_3(\XX) = \frac{1}{2}\sum_{i=1}^N\left( t_i - \frac{\sum_{j=1}^Nt_j}{N}\right)^2$,
\begin{itemize}
\item With this regularization term we aim at having lengths $t_i$, $1 \leq i \leq N$, to be close to their average length. This also forces the overall length to be equally distributed among segments.
\end{itemize}
\end{itemize}
Then we solve and compare
\begin{align}
\label{eq:len_all}
\min_{\XX}R_1(\XX) & \quad \mathrm{subject\ to}  \quad c_2(\XX) = 0\ , \\
\label{eq:pen_match_reg}
\min_{\XX}\left[ F_1(\XX) + R_1(\XX)\right] & \quad \mathrm{subject\ to}  \quad c_1(\XX) = 0\ , \\
\label{eq:gap_pen_reg}
\min_{\XX}\left[ F_2(\XX) + R_1(\XX)\right] & \quad \mathrm{subject\ to}  \quad c_3(\XX) = 0 \ , \\
\label{eq:gap_pen_reg_diff}
\min_{\XX}\left[ F_2(\XX) + R_2(\XX)\right] & \quad \mathrm{subject\ to}  \quad c_3(\XX) = 0 \ , \\
\label{eq:gap_pen_reg_var}
\min_{\XX}\left[ F_2(\XX) + R_3(\XX)\right] & \quad \mathrm{subject\ to}  \quad c_3(\XX) = 0 \ .
\end{align}
As in the previous case we will also consider the unconstrained version, therefore we solve
\begin{align}
\label{eq:unconst_reg}
\min_{\XX}\left[ F_3(\XX) + R_1(\XX)\right]\ .
\end{align}
Problem formulation \eqref{eq:len_all} is a regularized version of the constraint solving problem $c_2(\XX) = 0$. The remaining minimization problems are similar to problems in \eqref{eq:pen_match}-\eqref{eq:unconst}, however, we add a regularization term into the objective function $F(\XX)$.

We are interested in numerical behaviour in order to select the most suitable formulation of problem \eqref{eq:minF(X)} as minimization problems~\eqref{eq:len_all}-\eqref{eq:unconst_reg}. Note that there are many more variations than the ones shown. We use the sequential quadratic programming method~\cite[Ch.~18]{Nocedal:2006} to solve the minimization problems.

Some of the problem formulations will be removed from our considerations because of practical reasons. For example in Section \ref{sec:Dis_OFaC} it follows from Lemma \ref{lem:gxL} that \eqref{eq:gap_pen_reg} and \eqref{eq:unconst_reg} are not suitable.
\section{Review of Sequential Quadratic Programming}
\label{sec:SQP}
For reader's convenience we review the Sequential Quadratic Programming method. We are concerned with the following constrained minimization problem
\begin{equation}
\label{eq:arg_min}
\min_{x \in \Gamma} F(x)\ ,
\end{equation}
where a feasible set $\Gamma \subseteq \RR^n$ is defined by the system of equations
\begin{equation}
\label{eq:feasible_set}
\Gamma = \left\lbrace x \in \RR^n \mid c_k(x) = 0, 1 \leq k \leq m \right\rbrace,
\end{equation}
where $m \leq n$. We assume the functions $F: \RR^n \rightarrow \RR$ and $c_k: \RR^n \rightarrow \RR$, $1 \leq k \leq m$, to be twice continuously differentiable. We denote their gradients by $\nabla F(x)$, $\nabla c_k(x)$, $1 \leq k \leq m$, and their Hessian matrices by $\nabla^2 F(x)$, $\nabla^2 c_k(x)$, $1 \leq k \leq m$. For our convenience we use the vector notation $c(x) = \left[ c_1(x), \ldots, c_k(x) \right]^T \in \RR^m$, and for the Jacobian of constraints we put $B(x) = \left[ \nabla c_1(x), \ldots, \nabla c_m(x) \right] \in \RR^{n \times m}$. Let us suppose that the matrix $B(x)$ has full column rank.

We define the Lagrangian function as
\begin{equation}
\label{eq:Lagrangian}
\LL(x, \lambda) = F(x) + \lambda^T c(x)\ ,
\end{equation}
where $\lambda = \left[ \lambda_1, \ldots, \lambda_m \right]^T \in \RR^m$ is a vector of Lagrange multipliers.
The solution vector $x^\star \in \RR^n$ of \eqref{eq:arg_min} is said to satisfy the \emph{Karush-Kuhn-Tucker} (KKT) conditions, if and only if there exists a vector $\lambda^\star \in \RR^m$, such that
\begin{align}
\label{eq:gxLag}
\nabla_x \LL(x^\star, \lambda^\star) & = \nabla F(x^\star) + B(x^\star)\lambda^\star = 0\ , \\
\label{eq:glLag}
\nabla_\lambda \LL(x^\star, \lambda^\star) & = c(x^\star) = 0\ .
\end{align}
We denote the Hessian matrix of the Lagrangian by $\nabla_x^2 L(x, \lambda)$ and
\begin{equation}
\label{eq:hesLag}
\nabla_x^2 L(x, \lambda) = \nabla^2 F(x) + \sum_{k=1}^m \lambda_k \nabla^2c_k(x)\ .
\end{equation}
Then the \emph{second-order sufficient} conditions for a solution $x^\star \in \RR^n$ of \eqref{eq:arg_min} are
\begin{equation}
\label{eq:SSD}
w^T \nabla_x^2 L(x, \lambda) w \geq 0, \quad \mathrm{for\ all}\  w \in N\left(B^T\right)\ , 
\end{equation}
where $N\left(B^T\right) = \left\lbrace w \in \RR^n \mid B(x^\star)^Tw = 0 \right\rbrace$ and $\lambda^\star \in \RR^m$ satisfies the KKT conditions. For more details see \cite[Th. 12.5]{Nocedal:2006}.

We will use iterative methods for solving problem \eqref{eq:arg_min} and in each iteration we get
\begin{align}
\label{eq:x_update}
x_{\text{new}} = & \ x + \alpha_x d_x\ ,  \\
\label{eq:lam_update}
\lambda_{\text{new}} = & \ \lambda + \alpha_\lambda d_\lambda\ ,
\end{align}
where $d_x \in \RR^n$, $d_\lambda \in \RR^m$ are vectors, and $\alpha_x > 0$, $\alpha_\lambda > 0$ are step lengths. We use the Newton method to solve the KKT system of non-linear equations \eqref{eq:gxLag}-\eqref{eq:glLag} and get a system of  $n + m$ linear equations in $n + m$ unknowns, that is
\begin{align}
\label{eq:KKTsystem}
\begin{bmatrix}
H(x, \lambda) & B(x) \\
B(x)^T & 0
\end{bmatrix}
\begin{bmatrix}
d_x \\
d_\lambda
\end{bmatrix}
& =
\begin{bmatrix}
- \nabla_x \LL(x, \lambda) \\
- \nabla_\lambda \LL(x, \lambda) 
\end{bmatrix},
\end{align}
where $H(x, \lambda) \in \RR^{n\times n}$ is either $\nabla_{x}^2 \LL(x, \lambda)$ or an approximation of $\nabla_{x}^2 \LL(x, \lambda)$. We use the BFGS method, as described in \cite[p.~140]{Nocedal:2006}, for the approximation of the Hessian $\nabla_{x}^2 \LL(x, \lambda)$. Discussion and numerical experiments with using iterative methods to solve the KKT system  \eqref{eq:KKTsystem} can be found in \cite{Benzi:2005,LuksanVlcek:1998} and \cite{LuksanVlcek:2001}.
\section{Properties of Optimization Formulations}
\label{sec:Dis_OFaC}
In problem formulations \eqref{eq:pen_match}-\eqref{eq:unconst} we have no control over the lengths $t_i \geq t_0$, $ i = 1, \ldots, N$, of trajectories. Due to this, the computed error trajectory may feature degenerate segments, that is, a trajectory of zero length. It may even happen that during the algorithm we simulate the evolution backwards in time, that is, for some trajectories we  have $t_i < t_0$. 

This lack of control over the lengths causes numerical problems, especially, when we have many degenerate trajectories. In our experience, lengths of trajectories behave randomly. We do \emph{not} recommend to compute error trajectories using formulations \eqref{eq:pen_match}-\eqref{eq:unconst}.

To mitigate problems with random lengths of trajectories we introduced a regularization term in \eqref{eq:len_all}-\eqref{eq:unconst_reg}. Our goal is to distribute the lengths $t_i \geq t_0$, $i = 1, \ldots, N$, equally. However, we need to calculate the Lagrangian \eqref{eq:gxLag} for the solution vector $\XX \in \RR^{N(n+1)}$. When we compute the Lagrangian for \eqref{eq:len_all} -\eqref{eq:unconst_reg}, we observe hidden trouble. We address this problem in the following lemmata.

\begin{figure}[h]
\centering
\includegraphics[scale=1]{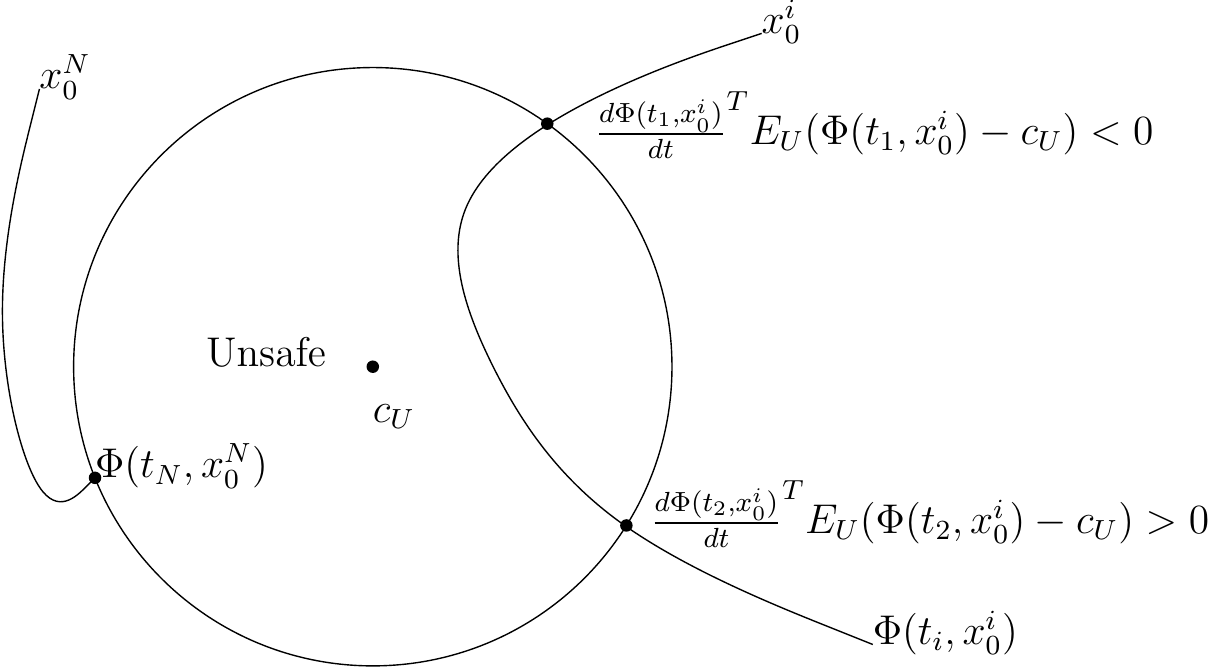}
\caption{An error trajectory enters and leaves $\Unsafe \subset \RR^n$.}
\label{fig:unsafe}
\end{figure}

First, we investigate the rank of the Jacobian of constraint. Lemma \ref{lem:Bx} concerns the Jacobian $B$ of constraints $c_1(\XX), c_2(\XX)$ and $c_3(\XX)$. Here we denote the sensitivity function $S: \RR\times\RR^n \to \RR^{n \times n}$ of the solution $x(t)= \left[ x_1(t), \ldots, x_n(t) \right]^T$ of the differential equation to the change of the initial value $x_0 = \left[x_{0,1}, \ldots, x_{0,n} \right]^T \in \RR^n$ by
\begin{equation}
\label{eq:Sensitivity}
S(t, x_0) = \frac{\partial \Phi(t, x_0)}{\partial x_0} = 
\begin{bmatrix}
 \frac{\partial x_1(t)}{\partial x_{0,1}} & \ldots & \frac{\partial x_1(t)}{\partial x_{0,n}} \\
  \frac{\partial x_2(t)}{\partial x_{0,1}} & \ldots & \frac{\partial x_2(t)}{\partial x_{0,n}} \\
\vdots & & \vdots \\
 \frac{\partial x_n(t)}{\partial x_{0,1}} & \ldots & \frac{\partial x_n(t)}{\partial x_{0,n}} \\
\end{bmatrix}\ .
\end{equation}
\begin{lemma}
\label{lem:Bx}
Let $\XX \in \RR^{N(n+1)}$ be a vector of parameters as in \eqref{eq:parameters}. 
Then
\begin{enumerate}
\item the Jacobian of the constraints $c_1(\XX)$, given by
\begin{equation}
\label{eq:grad_pen_match}
B = \begin{bmatrix}
-S(t_1, x_0^1)^T &  & & & \\
-\frac{d \Phi(t_1, x_0^1)}{d t_1}^T &   & & & \\
I &  - S(t_2, x_0^2)^T & & & \\
 & -\frac{d \Phi(t_2, x_0^2)}{d t_2}^T & & & \\
  &  I & \ddots & & \\
 &   & \ddots & & \\
 &  & & I & -S(t_{N-1}, x_0^{N-1})^T \\
 &  & &  & -\frac{d \Phi(t_{N-1}, x_0^{N-1})}{d t_{N-1}}^T \\
 &  & &  & I \\ 
  &  & &  & 0 \\ 
\end{bmatrix} \,,
\end{equation}
where $B \in \RR^{N(n+1) \times n(N-1)}$ has full column rank.
\item the Jacobian of  the constraints $c_3(\XX)$, given by
\begin{equation}
\label{eq:grad:gap_pen}
B = \begin{bmatrix}
E_I(x_0^1 - c_I) & 0 \\
0 & \vdots \\
0 & 0 \\
\vdots & S(t_N, x_0^N)^TE_U(\Phi(t_N, x_0^N) - c_U) \\
0 & \frac{d \Phi(t_{N}, x_0^{N})}{d t_{N}}^TE_U(\Phi(t_N, x_0^N) - c_U)
\end{bmatrix} \,,
\end{equation}
where $B \in \RR^{N(n+1) \times 2}$ has full column rank under the condition: $x_0^1 \neq c_I$ and there is at least one non-zero entry in the second column.
\item the Jacobian of the constraints $c_2(\XX)$, given by
\begin{equation}
\label{eq:grad_len_all}
B = \begin{bsmallmatrix}
E_I(x_0^1 - c_I) & -S(t_1, x_0^1)  & &   & & & \\
&-\frac{d \Phi(t_1, x_0^1)}{d t_1}^T &   &&  & & \\
&I &  - S(t_2, x_0^2)^T & & & & \\
 & & -\frac{d \Phi(t_2, x_0^2)}{d t_2}^T & & & & \\
 & &  I & \ddots & & & \\
& &   & \ddots & &  & \\
& &  & & I & -S(t_{N-1}, x_0^{N-1})^T  & \\
& &  & &  &  -\frac{d \Phi(t_{N-1}, x_0^{N-1})}{d t_{N-1}}^T & \\
& &  & &  & I &S(t_N, x_0^N)^TE_U(\Phi(t_N, x_0^N) - c_U)\\ 
 & &  & &  & 0 & \frac{d \Phi(t_{N}, x_0^{N})}{d t_{N}}^TE_U(\Phi(t_N, x_0^N) - c_U)
\end{bsmallmatrix} \,,
\end{equation}
where $B \in \RR^{N(n+1) \times n(N-1)+2}$ has full column rank under the condition: $x_0^1 \neq c_I$ and $\frac{d \Phi(t_{N}, x_0^{N})}{d t_{N}}^TE_U(\Phi(t_N, x_0^N) - c_U) \neq 0$.
\end{enumerate}
\begin{proof}
Proofs of the first and the third statement are in the Appendix in Lemmas~\ref{lem:grad_pen_match} and~\ref{lem:grad_len_all}. Full column rank of the Jacobian of constraints in \eqref{eq:grad:gap_pen} follows directly since matrix $E_I$ is symmetric positive definite.
\end{proof}
\end{lemma}
Let us discuss the conditions on the Jacobian in the second and third item of Lemma~\ref{lem:Bx}. The corresponding constraints ensure the initial state $x_0^1$ and the final state $\Phi(t_N, x_0^N)$ to be on the boundary of $\Init \subset \RR^n$, and $\Unsafe \subset \RR^n$ respectively. This implies $x_0^1 \neq c_I$, since $c_I$ is the centre of the set $\Init$. In order to fulfil the second part of the condition in the third item we need the term $\frac{d \Phi(t_{N}, x_0^{N})}{d t_{N}}^TE_U(\Phi(t_N, x_0^N) - c_U)$ to be non-zero. Whenever it does become zero during computation---which is unlikely---we change the size of the step $\alpha$ in \eqref{eq:x_update}-\eqref{eq:lam_update} to overcome that the Jacobian in \eqref{eq:grad_len_all} has linearly dependent columns. However, since an error trajectory enters the set $\Unsafe$ we expect the value to be negative, as shown in Figure \ref{fig:unsafe}. 

 Note that if the Jacobian does not have full column rank one can still solve the system in \eqref{eq:KKTsystem} if the right-hand side belongs to the range of the saddle point matrix \cite[p.~124]{Strang:2009}. 

\begin{figure}
\centering
\includegraphics[scale=1.3]{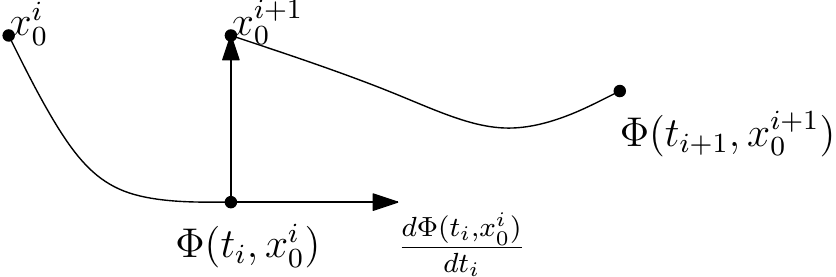}
\caption{A gap between two consecutive segments of an error trajectory.}
\label{fig:Gap}
\end{figure}

Lemma \ref{lem:gxL} concerns the Karush-Kuhn-Tucker conditions and the gradient of the Lagrangian
\[
\nabla_{\XX}\LL(\XX, \lambda) = \nabla F(\XX) + B\lambda\,.
\]
\begin{lemma}
\label{lem:gxL}
Depending on the problem formulation, the Lagrangian $\LL(\XX, \lambda)$ in \eqref{eq:gxLag} has the following form.
\begin{enumerate}
\item For problem formulation \eqref{eq:len_all} it is of the form
\begin{equation}
\label{eq:gxL_len_all}
\nabla_{\XX}\LL(\XX, \lambda)  = 
\begin{bmatrix}
\lambda_I E_I \left( x_0^1 - c_I\right) - S(t_1, x_0^1)^T \lambda_1 \\
t_1 - \frac{d\Phi(t_1, x_0^1)}{d t_1}^T\lambda_1 \\
\lambda_1 - S(t_2, x_0^2)^T\lambda_2 \\
t_2 - \frac{d\Phi(t_2, x_0^2)}{d t_2}^T\lambda_2 \\
\vdots \\
\lambda_{N-2} - S(t_{N-1}, x_0^{N-1})^T\lambda_{N-1} \\
t_{N-1} - \frac{d\Phi(t_{N-1}, x_0^{N-1}}{d t_{N-1}}^T\lambda_{N-1} \\
\lambda_{N-1} - \lambda_U S( t_N, x_0^N)^T E_U \left( \Phi(t_N, x_0^N) - c_U \right) \\
t_N + \lambda_U \frac{d \Phi(t_N, x_0^N)}{d t_N}^T E_U \left( \Phi(t_N, x_0^N) - c_U \right)
\end{bmatrix}\,,
\end{equation}
where $\lambda_i \in \RR^n$, $i = 1, \ldots, N-1$, and $\lambda_I \in \RR$, $\lambda_U \in \RR$.
\item For problem formulation \eqref{eq:pen_match_reg} it is of the form
\begin{equation}
\label{eq:gxL_pen_match_reg}
\nabla_{\XX}\LL(\XX, \lambda)  = 
\begin{bmatrix}
E_I \left( x_0^1 - c_I\right) - S(t_1, x_0^1)^T \lambda_1 \\
t_1 - \frac{d\Phi(t_1, x_0^1)}{d t_1}^T\lambda_1 \\
\lambda_1 - S(t_2, x_0^2)^T\lambda_2 \\
t_2 - \frac{d\Phi(t_2, x_0^2)}{d t_2}^T\lambda_2 \\
\vdots \\
\lambda_{N-2} - S(t_{N-1}, x_0^{N-1})^T\lambda_{N-1} \\
t_{N-1} - \frac{d\Phi(t_{N-1}, x_0^{N-1}}{d t_{N-1}}^T\lambda_{N-1} \\
\lambda_{N-1} -  S( t_N, x_0^N)^T E_U \left( \Phi(t_N, x_0^N) - c_U \right) \\
t_N +  \frac{d \Phi(t_N, x_0^N)}{d t_N}^T E_U \left( \Phi(t_N, x_0^N) - c_U \right)
\end{bmatrix}\,,
\end{equation}
where $\lambda_i \in \RR^n$, $i = 1, \ldots, N-1$.
\item For problem formulation \eqref{eq:gap_pen_reg} it is of the form
\begin{equation}
\label{eq:gxL_gap_pen_reg}
\nabla_{\XX}\LL(\XX, \lambda)  = 
\begin{bmatrix}
\lambda_I E_I \left( x_0^1 - c_I\right) - S(t_1, x_0^1)^T \left( x_0^2 - \Phi(t_1, x_0^1) \right) \\
t_1 - \frac{d\Phi(t_1, x_0^1)}{d t_1}^T\left( x_0^2 - \Phi(t_1, x_0^1) \right) \\
\left( x_0^2 - \Phi(t_1, x_0^1) \right) - S(t_2, x_0^2)^T\left( x_0^3 - \Phi(t_2, x_0^2) \right) \\
t_2 - \frac{d\Phi(t_2, x_0^2)}{d t_2}^T\left( x_0^3 - \Phi(t_2, x_0^2) \right) \\
\vdots \\
\left( x_0^{N-1} - \Phi(t_{N-2}, x_0^{N-2}) \right) - S(t_{N-1}, x_0^{N-1})^T\left( x_0^N - \Phi(t_{N-1}, x_0^{N-1}) \right) \\
t_{N-1} - \frac{d\Phi(t_{N-1}, x_0^{N-1}}{d t_{N-1}}^T\left( x_0^N - \Phi(t_{N-1}, x_0^{N-1}) \right) \\
\left( x_0^N - \Phi(t_{N-1}, x_0^{N-1}) \right) - \lambda_U S( t_N, x_0^N)^T E_U \left( \Phi(t_N, x_0^N) - c_U \right) \\
t_N + \lambda_U \frac{d \Phi(t_N, x_0^N)}{d t_N}^T E_U \left( \Phi(t_N, x_0^N) - c_U \right)
\end{bmatrix}\,,
\end{equation}
where $\lambda_I \in \RR$ and $\lambda_U \in \RR$.
\item For problem formulation \eqref{eq:gap_pen_reg_diff} it is of the form
\begin{equation}
\label{eq:gxL_gap_pen_reg_diff}
\nabla_{\XX}\LL(\XX, \lambda)  = 
\begin{bmatrix}
\lambda_I E_I \left( x_0^1 - c_I\right) - S(t_1, x_0^1)^T \left( x_0^2 - \Phi(t_1, x_0^1) \right) \\
-(t_2 - t_1) - \frac{d\Phi(t_1, x_0^1)}{d t_1}^T\left( x_0^2 - \Phi(t_1, x_0^1) \right) \\
\left( x_0^2 - \Phi(t_1, x_0^1) \right) - S(t_2, x_0^2)^T\left( x_0^3 - \Phi(t_2, x_0^2) \right) \\
(t_2 - t_1) - (t_3 - t_2) - \frac{d\Phi(t_2, x_0^2)}{d t_2}^T\left( x_0^3 - \Phi(t_2, x_0^2) \right) \\
\vdots \\
\left( x_0^{N-1} - \Phi(t_{N-2}, x_0^{N-2}) \right) - S(t_{N-1}, x_0^{N-1})^T\left( x_0^N - \Phi(t_{N-1}, x_0^{N-1}) \right) \\
(t_{N-1} - t_{N-2}) - (t_N - t_{N-1}) - \frac{d\Phi(t_{N-1}, x_0^{N-1}}{d t_{N-1}}^T\left( x_0^N - \Phi(t_{N-1}, x_0^{N-1}) \right) \\
\left( x_0^N - \Phi(t_{N-1}, x_0^{N-1}) \right) - \lambda_U S( t_N, x_0^N)^T E_U \left( \Phi(t_N, x_0^N) - c_U \right) \\
(t_N - t_{N-1}) + \lambda_U \frac{d \Phi(t_N, x_0^N)}{d t_N}^T E_U \left( \Phi(t_N, x_0^N) - c_U \right)
\end{bmatrix}\,,
\end{equation}
where $\lambda_I \in \RR$ and $\lambda_U \in \RR$.
\item For problem formulation \eqref{eq:gap_pen_reg_var} it is of the form
\begin{equation}
\label{eq:gxL_gap_pen_reg_var}
\nabla_{\XX}\LL(\XX, \lambda)  = 
\begin{bmatrix}
\lambda_I E_I \left( x_0^1 - c_I\right) - S(t_1, x_0^1)^T \left( x_0^2 - \Phi(t_1, x_0^1) \right) \\
t_1 - \bar{t} - \frac{d\Phi(t_1, x_0^1)}{d t_1}^T\left( x_0^2 - \Phi(t_1, x_0^1) \right) \\
\left( x_0^2 - \Phi(t_1, x_0^1) \right) - S(t_2, x_0^2)^T\left( x_0^3 - \Phi(t_2, x_0^2) \right) \\
t_2 - \bar{t} - \frac{d\Phi(t_2, x_0^2)}{d t_2}^T\left( x_0^3 - \Phi(t_2, x_0^2) \right) \\
\vdots \\
\left( x_0^{N-1} - \Phi(t_{N-2}, x_0^{N-2}) \right) - S(t_{N-1}, x_0^{N-1})^T\left( x_0^N - \Phi(t_{N-1}, x_0^{N-1}) \right) \\
t_{N-1}- \bar{t} - \frac{d\Phi(t_{N-1}, x_0^{N-1}}{d t_{N-1}}^T\left( x_0^N - \Phi(t_{N-1}, x_0^{N-1}) \right) \\
\left( x_0^N - \Phi(t_{N-1}, x_0^{N-1}) \right) - \lambda_U S( t_N, x_0^N)^T E_U \left( \Phi(t_N, x_0^N) - c_U \right) \\
t_N - \bar{t} + \lambda_U \frac{d \Phi(t_N, x_0^N)}{d t_N}^T E_U \left( \Phi(t_N, x_0^N) - c_U \right)
\end{bmatrix}\,,
\end{equation}
where $\lambda_I \in \RR$, $\lambda_U \in \RR$ and $\bar{t} = \frac{1}{N}\sum_{i=1}^N t_i$.
\item For problem formulation \eqref{eq:unconst_reg} it is of the form
\begin{equation}
\label{eq:gxL_unconst_reg}
\nabla_{\XX}\LL(\XX, \lambda)  = 
\begin{bmatrix}
E_I \left( x_0^1 - c_I\right) - S(t_1, x_0^1)^T \left( x_0^2 - \Phi(t_1, x_0^1) \right) \\
t_1 - \frac{d\Phi(t_1, x_0^1)}{d t_1}^T\left( x_0^2 - \Phi(t_1, x_0^1) \right) \\
\left( x_0^2 - \Phi(t_1, x_0^1) \right) - S(t_2, x_0^2)^T\left( x_0^3 - \Phi(t_2, x_0^2) \right) \\
t_2 - \frac{d\Phi(t_2, x_0^2)}{d t_2}^T\left( x_0^3 - \Phi(t_2, x_0^2) \right) \\
\vdots \\
\left( x_0^{N-1} - \Phi(t_{N-2}, x_0^{N-2}) \right) - S(t_{N-1}, x_0^{N-1})^T\left( x_0^N - \Phi(t_{N-1}, x_0^{N-1}) \right) \\
t_{N-1} - \frac{d\Phi(t_{N-1}, x_0^{N-1}}{d t_{N-1}}^T\left( x_0^N - \Phi(t_{N-1}, x_0^{N-1}) \right) \\
\left( x_0^N - \Phi(t_{N-1}, x_0^{N-1}) \right) - S( t_N, x_0^N)^T E_U \left( \Phi(t_N, x_0^N) - c_U \right) \\
t_N + \frac{d \Phi(t_N, x_0^N)}{d t_N}^T E_U \left( \Phi(t_N, x_0^N) - c_U \right)
\end{bmatrix}\,.
\end{equation}
\end{enumerate}
\begin{proof}
We obtain these results directly after substituing into the formula $\nabla_{\XX}\LL(\XX, \lambda) = \nabla F(\XX) + B\lambda$. Here we take $B$ from Lemma \ref{lem:Bx}. Depending on the vector of constraints we get $\lambda$ that is either $[\lambda_I, \lambda_U]\in \RR^2$, or $[\lambda_1, \ldots, \lambda_{N-1}]\in \RR^{n(N-1)}$, or $[\lambda_I,\lambda_1, \ldots, \lambda_{N-1}, \lambda_U]\in \RR^{n(N-1) + 2}$. In the unconstrained case the term $B\lambda$ vanishes.
\end{proof}
\end{lemma}
From Lemma \ref{lem:gxL} it follows that introducing the regularization term $\sum_{i=1}^N t_i^2$ into the objective function may prevent obtaining an error trajectory: 
As illustrated in Figure \ref{fig:Gap}, when the distance $\|  x_0^{i+1} - \Phi(t_i, x_0^i) \|_2^2$ is minimal with respect to time $t_i$, then the vectors
\[
\frac{d\Phi(t_i, x_0^i)}{d t_i} \quad \text{and} \quad  x_0^{i+1} - \Phi(t_i, x_0^i)
\]
are perpendicular for $1 \leq i \leq N-1$. In this case, the terms $\frac{d\Phi(t_i, x_0^i)}{d t_i}^T \left( x_0^{i+1} - \Phi(t_i, x_0^i) \right)$ in the third and sixth item of Lemma \ref{lem:gxL} are zero. This forces the lengths $t_i$ to be zero.  Moreover, also the final goal of fulfilling the matching conditions  $x_0^{i+1} = \Phi(t_i, x_0^i)$, $1 \leq i \leq N-1$ has the same effect, forcing the lengths $t_i$ to be zero. Because of these reasons we do not recommend using \eqref{eq:gap_pen_reg} and \eqref{eq:unconst_reg} for computing error trajectories. 
\section{Practical Considerations}
\label{sec:Numerics}
When we select the problem formulation we also consider the structure and sparsity of the saddle point matrix in \eqref{eq:KKTsystem}. The form of $B$, that is, the Jacobian of the constraints vector, is described in Lemma \ref{lem:Bx}. We can influence the structure of the saddle-point matrix by choosing an approximation scheme for the Hessian of the Lagrangian. We aim at keeping its sparse structure, that is, we want to avoid having the Hessian to be a dense matrix.
\subsection{Hessian of the Lagrangian}
\label{subsec:HoL}
It is convenient to use the BFGS method~\cite[p.~140]{Nocedal:2006}. When we set $s = \XX_{new} - \XX$ and $y = \nabla_{\XX}\LL(\XX_{new},\lambda_{new}) - \nabla_{\XX}\LL(\XX, \lambda_{new})$, where $\XX, \XX_{new}, \lambda$ and $\lambda_{new}$ is from \eqref{eq:x_update}, \eqref{eq:lam_update}, then the BFGS updating scheme is given by
\begin{equation}
\label{eq:H_BFGSupdate}
H_{new} = H - \frac{Hss^TH}{s^THs} + \frac{yy^T}{y^Ts}\,,
\end{equation}
where the initial approximation of the Hessian is the identity matrix in our implementation. However, in this way we obtain a dense approximation $H(\XX, \lambda) \in \RR^{N(n+1) \times N(n+1)}$.

When we solve the minimization problem~\eqref{eq:len_all}, then the Jacobian of constraints $B$ is given in \eqref{eq:grad_len_all} and the Hessian is a block-diagonal matrix such that
\[
H = \begin{bmatrix}
H_1 & &   \\
   & \ddots &  \\
   & & H_N 
\end{bmatrix}\,,
\]
where $N$ is the number of trajectories, $H_i \in \RR^{(n+1) \times (n+1)}$, $1 \leq i \leq N$. Therefore we may use the BFGS method block by block, keeping the block-diagonal structure \cite{Griewank:1982}. In this fashion we obtain a block-diagonal symmetric definite approximation of the Hessian. Moreover, the matrix is sparse and the ratio of non-zero elements to zeros is $1/N$. The block-diagonal structure also appears when we consider problem formulation \eqref{eq:pen_match_reg} with $B$ in \eqref{eq:grad_pen_match}.

When we use a different regularization term such as $\frac{1}{2}\sum_{i=1}^{N-1}(t_{i+1} - t_i)^2$ as in \eqref{eq:gap_pen_reg_diff}, then the parameters are no longer separable, however, they are partially separable and the resulting Hessian is a banded matrix. If we apply the BFGS updating scheme from \eqref{eq:H_BFGSupdate}, then we obtain a dense approximation of the Hessian matrix $H(\XX, \lambda)$. However, we may again use the BFGS method block-wise~\cite{Griewank:1982}, although, this time we work with blocks of size $2(n+1) \times 2(n+1)$. As a result we obtain a banded symmetric positive definite approximation of the Hessian. Its ratio of non-zero elements to zeros is approximately $4/N$.

Another regularization term we consider is $\frac{1}{2}\sum_{i=1}^N\left( t_i - \sum_{j=1}^Nt_j/N\right)^2$, however, this effectively connects all parameters. Therefore, it leaves us only with a dense approximation when using BFGS. 

\subsection{Numerical Solution of the KKT system}
\label{subsec:KKT}
We also need to address  solution techniques for the KKT system in \eqref{eq:KKTsystem}. There are several approaches. Two of them are the \emph{Schur-complement} method and the \emph{Null-space} method. For a thorough overview of methods for solving the KKT system, also known as a \emph{saddle point} problem, see~\cite[p.~29-59]{Benzi:2005}.

Note that it may happen when using structure preserving approximation for the Hessian, that its condition number gets worse than using \eqref{eq:H_BFGSupdate} instead. In our experience the approximation of the Hessian is ill-conditioned when we use the BFGS method from \eqref{eq:H_BFGSupdate}. Because of that we cannot use the Schur-complement method. We can apply the Null-space method which works for a singular Hessian approximation as long as its projection on the null-space $N(B^T)$ of $B^T$ is symmetric positive definite. This is the second order sufficient condition on $[\XX^T, \lambda^T ]^T$ to be the solution of \eqref{eq:KKTsystem}. 

Denote by $N_{B^T}$ a matrix whose columns form an orthonormal basis of $N(B^T)$. During numerical testing we noticed that the condition number of the projected Hessian $N_{B^T}^THN_{B^T}$ tends to be the same as the condition number of the Hessian $H$ in problem formulations \eqref{eq:gap_pen_reg_diff} and \eqref{eq:gap_pen_reg_var}, when we project the Hessian on the null-space of $B^T$ in \eqref{eq:grad:gap_pen}. When we project the Hessian on the null-space of $B^T$ from \eqref{eq:grad_len_all} and \eqref{eq:grad_pen_match} its condition number is usually several magnitudes lower than the condition number of the Hessian. 

Since the Jacobian $B$ is a sparse matrix one may use sparse QR decomposition to compute $N_{B^T}$ as proposed in \cite[Alg.~1]{Foster:1986}. Another possibility is to use Givens rotations on non-zero elements of $B$ to compute its QR-decomposition \cite[p.~227]{Golub:1996}. One can also avoid the computation of the basis $N(B^T)$ at all by \cite{Gould:2001}, \cite[Alg.~5.1.3]{Dollar:2005} and the preconditioned projected conjugate gradient method \cite[Alg.~NPCG]{LuksanVlcek:2001}.

\begin{figure}
\centering
\includegraphics[scale=0.25]{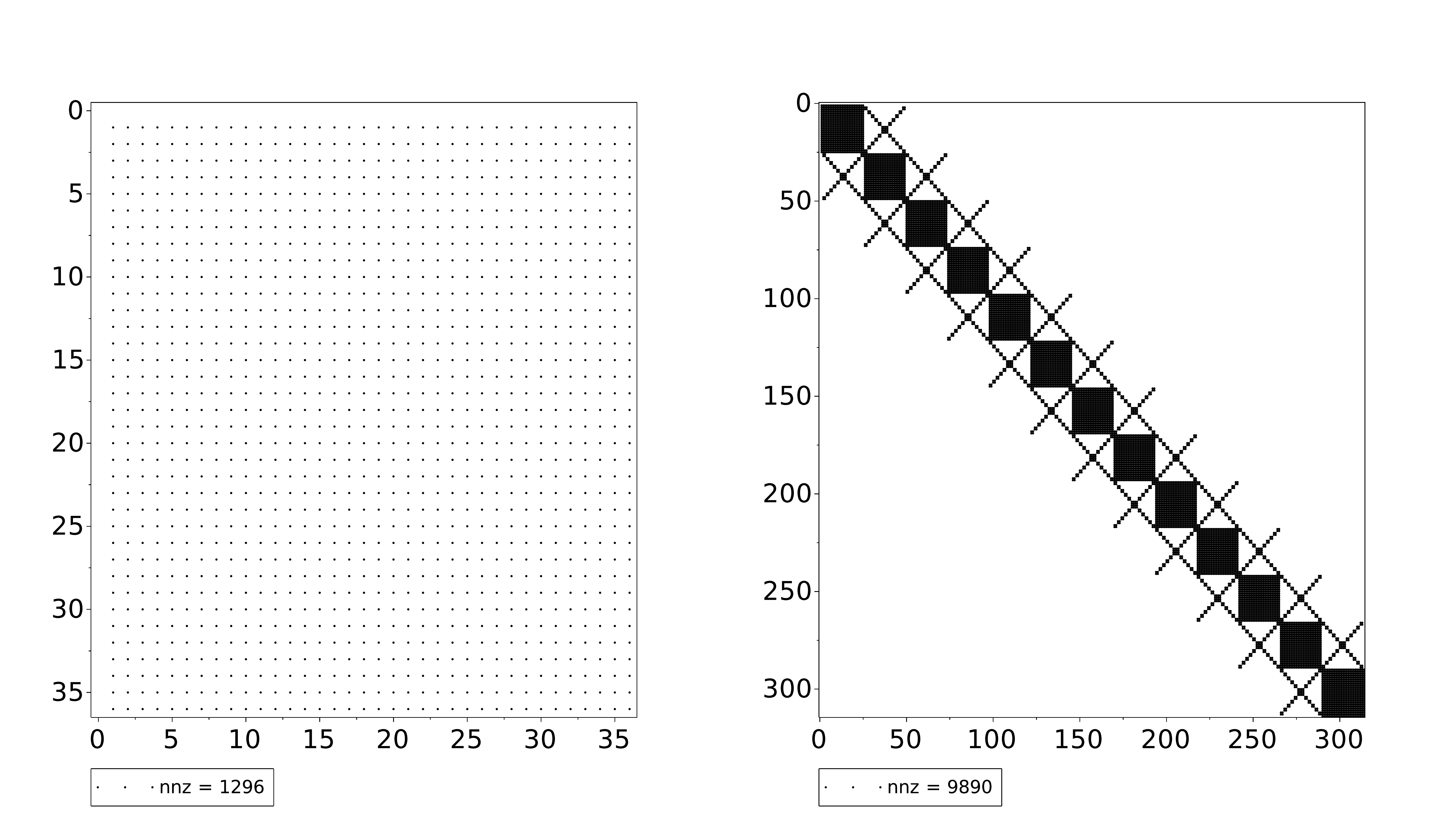}
\caption{Visualization of matrices $N_{B^T}^THN_{B^T}$ on the left and $B^TB$ on the right. Data come from Benchmark \ref{subsec:B1}, problem formulation \eqref{eq:len_all} for $n = 24$ and $N = 14$ with the block-diagonal Hessian matrix.}
\label{fig:Struc_Bx'Bx_N'HN}
\end{figure}

Whether we apply QR-decomposition to compute the orthonormal basis of $N(B^T)$ or avoid its computation depends on the size and the structure of matrices $B^TB$ and $N_{B^T}^THN_{B^T}$. Note that from QR-decomposition \cite[p.~227]{Golub:1996} one gets dense basis vectors of $N(B^T)$ for Jacobian matrices \eqref{eq:grad_pen_match} and \eqref{eq:grad_len_all}. Therefore the matrix $N_{B^T}^THN_{B^T}$ is symmetric positive definite and dense. However, the matrix $B^TB$ is symmetric positive definite and banded. One can see these structures in Figure \ref{fig:Struc_Bx'Bx_N'HN}. From our experience, when the number of segments $N$ is much larger than the state-space dimension $n$ it is preferable to avoid computation of $N(B^T)$.
\section{Computational Experiments}
\label{sec:Benchmarks}
In this section we describe the  algorithm that we use for finding error trajectories of ordinary differential equations. There are two steps. First, we find a candidate for an error trajectory by using SQP as shown in Algorithm \ref{alg:fals_SQP}. Second, we verify the result by simulation as described in Section \ref{subsec:Methodology}. In addition we present a series of benchmarks where we compare several different problem formulations. For our computation we use \emph{Scilab 5.5.2} running on \emph{Cent OS 6.8}. 

\begin{algorithm}[h]
\label{alg:fals_SQP}
\KwData{Problem formulation for minimization problem in \eqref{eq:minF(X)}, and initial vector of parameters $\XX \in \RR^{N(n+1)}$}
\KwResult{A new vector $\XX$ that corresponds to an error trajectory candidate}
\While{stopping criteria are not met}
{
Approximate the Hessian $H(\XX, \lambda)$\;
Compute the Jacobian of constraints and construct the KKT system\;
Solve the KKT system for $d_{\XX}$ and $d_{\lambda}$ \;
Select a step-size $\alpha$ using a merit function\;
Set $\XX = \XX + \alpha d_{\XX}$, $\lambda = \lambda + \alpha d_{\lambda}$\;
}
\caption{SQP method for finding error trajectories}
\end{algorithm}

In Algorithm \ref{alg:fals_SQP} we use the BFGS method as described in \eqref{eq:H_BFGSupdate} to approximate $H(\XX, \lambda)$. Whenever the descent condition $y^Ts > 0$ is not satisfied we use the current approximation of the Hessian without any update. Another possibility to consider is applying \emph{damped} BFGS \cite[p.~537]{Nocedal:2006} to avoid skipping the update of the Hessian. In our experience this does not lead to an improvement. To compute the Jacobian of constraints we solve \emph{variational equations} \cite[Sec.~2.4]{Ascher:1998} in order to get the sensitivity functions \eqref{eq:Sensitivity}. For solving the KKT system for $d_\XX$ and $d_\lambda$ we choose the preconditioned projected conjugate gradient method \cite[Alg.~NPCG]{LuksanVlcek:2001} that avoids computation of the null-space basis of $B^T$. In our implementation of \cite[Alg.~NPCG]{LuksanVlcek:2001} we set matrix $D$ to be the identity matrix, therefore, the constraint preconditioner we use is of the form
\[
C
=
\begin{bmatrix}
I & B(x) \\
B(x)^T & 0
\end{bmatrix}\,.
\]
We select a step-size $\alpha = \alpha_x = \alpha_\lambda$ by the line search method, making the value of the following merit function decrease in each iteration
\begin{equation}
\label{eq:merit}
m(\alpha)  = F(\XX + \alpha d_\XX) + (\lambda + d_\lambda)^Tc(\XX + \alpha d_\XX) + \frac{\omega}{2}\|c(\XX + \alpha d_x)\|_2^2\ ,
\end{equation}
where $\omega$ is a parameter. Properties of the merit function in \eqref{eq:merit} are formulated in \cite[Th.~2.5]{LuksanVlcek:1998} and  \cite[Th.~8]{LuksanVlcek:2001}. 
The step-size $\alpha$ is accepted for $\XX_{new} = \XX + \alpha d_\XX$ and $\lambda_{new} = \lambda + \alpha d_\lambda$ when 
\begin{equation}
\label{eq:stepsize_accept}
m(\alpha) - m(0) \leq \delta \alpha m'(0)\ ,
\end{equation} 
where $\delta$ is a parameter, and 
the derivative $m'(0)$ of the merit function at $\alpha = 0$ is given by 
\begin{equation}
\label{eq:merit_der}
m'(0) = d_\XX^T\left( \nabla F(\XX) + B(\lambda + d_\lambda) \right) + \omega d_\XX^TBc(\XX)\ .
\end{equation}

For all benchmarks we set $\omega = 1$ and $\delta = 1\times 10^{-4}$. In the end, we use the built in function \emph{ode} \cite{Scilab} for solving differential equations. The \emph{ode} solver in default setting calls the \emph{lsoda} solver of the package ODEPACK. It automatically selects between stiff and non-stiff methods.

We terminate Algorithm \ref{alg:fals_SQP} whenever one of these stopping criteria is met:
\begin{enumerate}
\item[] S 1: $\| \nabla_\XX \LL(\XX, \lambda) \|_2 < \varepsilon_1$ and $\| c(\XX) \|_2 < \varepsilon_2$,
\item[] S 2: the maximum number of iteration MAXIT is reached,
\item[] S 3: the step-size $\alpha < \varepsilon_3$,
\end{enumerate}
where we put $\varepsilon_1 = 1\times 10^{-3}$, $\varepsilon_2 = 1\times 10^{-8}$, $\varepsilon_3 = 1\times 10^{-8}$ and MAXIT to be $400$ iterations.

\subsection{Methodology}
\label{subsec:Methodology}
In a series of benchmarks we undertake the following procedure. We choose a point $c_I \in \RR^n$ and solve the system of differential equations in \eqref{eq:DiffEq} with $c_I$ to be the initial condition. We consider the time interval $0 \leq t \leq 5$. We denote the end-state of a computed trajectory by $c_U \in \RR^n$. Initial and unsafe sets of states are then $n$-dimensional balls centred at $c_I$, and $c_U$ respectively. The radius of these balls is equal to $1/4$ each.

Once we have created $\Init$ and $\Unsafe$, we proceed with splitting our trajectory into $N$ segments of the same lengths. We mark initial and final states by $x_0^i$ and $\Phi(t_i, x_0^i)$ with $t_i = 5/N$ for $1 \leq i \leq N$. Then we change initial states according to the following rule. For all $1 \leq i \leq N$ we update initial states $x_0^i + u$, where  $u = 0.5\times\left[-1, 1, \ldots, (-1)^n\right]^T \in \RR^n$. With these updated initial conditions and lengths we form a new vector of parameters $\XX \in \RR^{N(n+1)}$.

We run Algorithm \ref{alg:fals_SQP} and obtain a new vector of parameters that corresponds to an error trajectory candidate consisting of $N$ segments.  In order to verify our result we simulate the system for $t = \sum_{i=1}^N t_i$ time units originating in $x_0^1$. We check whether for this newly computed states $\| x_0^1 - c_I \|_{E_I} \leq 1 + \varepsilon_4$ and $\| \Phi(t,x_0^1) - c_U \|_{E_U} \leq 1 + \varepsilon_4$ hold, where $\varepsilon_4 = 1\times 10^{-4}$. If these two inequalities are satisfied we call such a trajectory an error trajectory and our method succeeded. If Scilab fails to solve the ODE or our computed trajectory is not an error trajectory we mark the corresponding row in the tables with results with the flag ``F''. Since we do not put any restriction on the lengths $t_i$, $1 \leq i \leq N$, and since the \emph{ODE} solver in Scilab is able to simulate the evolution backward in time we also mark by ``F'' those results that have at least one segment with negative length.

Let us discuss the choice of problem formulation for the following benchmark problems. We aim at comparing various choices for objective function and constraints as well as different approximation schemes for the Hessian matrix. To this end we chose problem formulations \eqref{eq:len_all}, \eqref{eq:pen_match_reg} and \eqref{eq:gap_pen_reg_diff}. Since the Hessian matrix of the Lagrangian for  \eqref{eq:len_all} and \eqref{eq:pen_match_reg} is block-diagonal, we compare full BFGS approximation with BFGS approximation applied block-wise. The Hessian matrix for \eqref{eq:gap_pen_reg_diff} is banded, hence, we compare full BFGS approximation with BFGS approximation that keeps the banded structure. In the tables that follow we put results corresponding to the structured approximation of the Hessian matrix in the left, and those corresponding to full BFGS approximation in the right respectively.

Finally, we do not try problem formulations \eqref{eq:gap_pen_reg} and \eqref{eq:unconst_reg} because from Lemma \ref{lem:gxL} it follows that result is usually not an error trajectory. Moreover, we omit problem formulation  \eqref{eq:gap_pen_reg_var} since preliminary tests showed that the results are similar to problem formulation \eqref{eq:gap_pen_reg_diff}.
\subsection{Benchmark 1}
\label{subsec:B1}
We consider the following benchmark problem where the dynamics is given by
\begin{equation}
\label{eq:B1}
\dot{x}(t) = Ax(t) + \sin(x^r(t)),
\end{equation}
where $x^r = [x_n(t), x_{n-1}(t), \ldots, x_1(t)]^T \in \RR^n$ and matrix $A$ is a block diagonal matrix such that
\[
A =
\begin{bmatrix}
0 & 1 & & &\\
-1 & 0 & & & \\
& & \ddots & & \\
& & & 0 & 1\\
& & & -1 & 0 
\end{bmatrix} \in \RR^{n \times n}\;.
\]
We use a matrix-vector notation, therefore, we read $\sin(x^r) = [\sin(x_n), \ldots \sin(x_1)]^T \in \RR^n$ and put $c_I = [1, \ldots, 1]^T \in \RR^n$.
\begin{table}
\centering
\begin{tabular}{*{24}{c}}
 n & N & NIT & S & n & N & NIT & S & n & N & NIT & S & n & N & NIT & S \\
\hline
 $10$  &  $5$  &  $37$  &  $1$  &  $20$  &  $5$  &  $54$  &  $1$  &  $30$  &  $5$  &  $46$  & $1$  &  $40$  &  $5$  &  $67$  &  $1$ \\
  &  $10$  &  $47$  &  $1$  &   &  $10$  &  $54$  &  $3$  &   &  $10$  &  $122$  &  $1$  &  &  $10$  &  $57$  &  $1$ \\
  &  $15$  &  $126$  &  $1$  &   &  $15$  &  $52$  &  $1$  &   &  $15$ &  $97$  &  $1$  &  &  $15$  &  $60$  &  $1$ \\
  &  $20$  &  $79$  &  $3$  & &  $20$  &  $149$  &  $3$  &   &  $20$  &  $101$  &  $3$  &  &  $20$ &  $90$  &  $3$ \\
  &  $25$  &  $185$  &  $1$  &   &  $25$  &  $107$  &  $3$ &   &  $25$  &  $80$  &  $1$  &  &  $25$  &  $98$  &  $3$ \\
  &  $30$  & $400$  &  $2$  &   &  $30$  &  $108$  &  $1$  &   &  $30$  &  $122$  &  $1$  &  &  $30$  &  $98$  &  $1$ \\[5mm]
n & N & NIT & S & n & N & NIT & S & n & N & NIT & S & n & N & NIT & S \\
\hline
 $10$  &  $5$  &  $33$  &  $1$  &  $20$  &  $5$  &  $38$  &  $1$  &  $30$  &  $5$  &  $38$ &  $1$  &  $40$  &  $5$  &  $39$  &  $1$ \\
   &  $10$  &  $43$  &  $1$  &   & $10$  &  $59$  &  $3$  &  &  $10$  &  $59$  &  $1$  & &  $10$  &  $70$ &  $1$ \\
   &  $15$  &  $66$  &  $3$  &   &  $15$  &  $67$  &  $3$  &  & $15$  &  $88$  &  $1$  & &  $15$  &  $86$  &  $3$ \\
   &  $20$  &  $75$  & $3$  &   &  $20$  &  $83$  &  $3$  &  &  $20$  &  $100$  &  $3$  &&  $20$  &  $103$  &  $3$ \\
   &  $25$  &  $80$  &  $1$  &   &  $25$  &  $98$ &  $1$  &  &  $25$  &  $115$  &  $1$  & &  $25$  &  $126$  &  $1$ \\
   &  $30$  &  $96$  &  $3$  &   &  $30$  &  $113$  &  $1$  &  &  $30$  & $126$  &  $3$  &    &  $30$  &  $137$  &  $3$ \\
\end{tabular}
\caption{Benchmark 1 problem formulation \eqref{eq:len_all}: there are results for block diagonal BFGS approximation in the top and the full BFGS approximation in the bottom.}
\label{tab:B1PF03BFGS}
\end{table}
\begin{table}
\centering
\begin{tabular}{*{24}{c}}
 n & N & NIT & S & n & N & NIT & S & n & N & NIT & S & n & N & NIT & S \\
\hline
 $10$  &  $5$  &  $57$  &  $1$  &  $20$  &  $5$  &  $170$  &  $1$  &  $30$  &  $5$  &  $108$ &  $1$  &  $40$  &  $5$  &  $223$  &  $1$ \\
   &  $10$  &  $61$  &  $1$  &  & $10$  &  $89$  &  $1$  &  &  $10$  &  $98$  &  $1$  &    &  $10$  &  $154$  & $1$ \\
   &  $15$  &  $79$  &  $1$  &  &  $15$  &  $400$  &  $2$  &  &  $15$  &  $106$  &  $1$  &  &  $15$  &  $114$  &  $1$ \\
   &  $20$  &  $117$  & $1$  &  &  $20$  &  $106$  &  $1$  &  &  $20$  &  $100$  &  $3$  &  & $20$  &  $108$  &  $3$ \\
   &  $25$  &  $85$  &  $1$  &  &  $25$  &  $121$  & $1$  &  &  $25$  &  $122$  &  $1$  &  &  $25$  &  $124$  &  $1$ \\
   &  $30$  &  $107$  &  $1$  &    &  $30$  &  $119$  &  $1$  &    &  $30$  &  $131$ &  $1$  &  &  $30$  &  $129$  &  $1$ \\[5mm]
n & N & NIT & S & n & N & NIT & S & n & N & NIT & S & n & N & NIT & S \\
\hline
 $10$  &  $5$  &  $86$  &  $1$  &  $20$  &  $5$  &  $400$  &  $2$  &  $30$  &  $5$  &  $400$ &  $2$  &  $40$  &  $5$  &  $75$  &  $1$ \\
   &  $10$  &  $96$  &  $1$  &   &  $10$  &  $400$  &  $2$  &    &  $10$  &  $108$  &  $1$  &    &  $10$  &  $400$  &  $2$ \\
   &  $15$  &  $121$  &  $1$  &    &  $15$  &  $171$  &  $1$  &   &  $15$  &  $154$  &  $1$  &   &  $15$  &  $166$  &  $1$ \\
   &  $20$  & $88$  &  $1$  &    &  $20$  &  $108$  &  $1$  &    &  $20$  &  $133$  &  $3$ &    &  $20$  &  $159$  &  $1$ \\
   &  $25$  &  $103$  &  $1$  &    &  $25$  &  $128$  &  $1$  &    &  $25$  &  $193$  &  $1$  &   &  $25$  &  $155$  &  $3$ \\
   &  $30$  &  $117$  &  $1$  &    &  $30$  &  $226$  &  $1$  &    &  $30$  &  $153$  &  $1$  &    &  $30$  &  $258$  &  $1$ \\ 
\end{tabular}
\caption{Benchmark 1 problem formulation \eqref{eq:pen_match_reg}: there are results for block diagonal BFGS approximation in the top and the full BFGS approximation in the bottom. }
\label{tab:B1PF04BFGS}
\end{table}
\begin{table}
\centering
\begin{tabular}{*{24}{c}}
 n & N & NIT & S & n & N & NIT & S & n & N & NIT & S & n & N & NIT & S \\
\hline
 $10$  &  $5$  &  $400$  &   F       &  $20$  &  $5$  &  -  &  F                       &  $30$  &  $5$  &  $400$  &  F  &  $40$  &  $5$  &  $400$  &  F \\
  &  $10$  &  $400$  &  $2$  &    &  $10$ &  $400$  &  F      &   &  $10$  &  -  &  F          &    &  $10$  &  -  &  F \\
  & $15$  &  $400$  &  F      &    &  $15$  &  $400$  &  F        &    &  $15$  &  -  & F            &    &  $15$  &  -  &  F \\
  &  $20$  &  $400$  &  $2$  &  &  $20$  & $400$  &  F       &   &  $20$  &  -  &  F          &    &  $20$  &  -  &  F \\
  & $25$  &  $400$  &  $2$  &    &  $25$  &  $400$  &  $2$  &    &  $25$  &  $400$ &  F   &    &  $25$  &  $400$  &  F \\
  &  $30$  &  $400$  &  $2$  &    & $30$  &  $400$  &  F    &   &  $30$  &  $400$  &  F    &    &  $30$  &  $400$ &  F \\ [5mm]
n & N & NIT & S & n & N & NIT & S & n & N & NIT & S & n & N & NIT & S \\
\hline
 $10$  &  $5$  &  $64$  &  $1$  &  $20$  &  $5$  &  $215$  &  F        &  $30$  &  $5$  &  $332$ &  F      &  $40$  &  $5$  &  -  &  F \\
  &  $10$  &  $182$  &  F      &    &  $10$ &  $400$  &  F    &    &  $10$  &  -  &          F      &    &  $10$  &  -  &  F \\
 &  $15$  &  $400$  &  F       &    &  $15$  &  -  &  F        &    &  $15$  &  -  &              F       &    &  $15$  &  -  &  F \\
  &  $20$  &  $266$  &  $1$  &    &  $20$  &  $400$  &  F      &    &  $20$  &  $400$  &  F  &   &  $20$  &  -  &  F \\
  & $25$  &  $360$  &  $1$  &    &  $25$  &  $358$  &  F       &    &  $25$  &  $400$ &  F  &    &  $25$  &  $400$  &  F \\
  &  $30$  &  $341$  &  F     &   &  $30$  &  $400$  &  F       &    &  $30$  &  $394$  &  F  &    &  $30$  &  $400$  &  $2$ \\
\end{tabular}
\caption{Benchmark 1 problem formulation \eqref{eq:gap_pen_reg_diff} there are results for block diagonal BFGS approximation in the top and the full BFGS approximation in the bottom. }
\label{tab:B1PF09BFGS}
\end{table}

For problem formulation \eqref{eq:len_all} there are results in Tab. \ref{tab:B1PF03BFGS}. One can see that the desired solution was computed every time. Note that the block-diagonal BFGS outperforms standard BFGS in the number of iterations only for $n = 40$. Also note that with the increasing number of segments our method needs more iteration when standard BFGS scheme for a dense approximation of the Hessian is used.

In Tab. \ref{tab:B1PF04BFGS} there are results for problem formulation \eqref{eq:pen_match_reg}. In this case we were able to compute the desired solution every time, however, the number of iterations required was higher than for \eqref{eq:len_all}. 

To conclude this part we show problem formulation \eqref{eq:gap_pen_reg_diff} and its results in Tab. \ref{tab:B1PF09BFGS}. There are many failed attempts marked by ``F''. Also, the \emph{dash} symbol shows when the \emph{ode} solver in Scilab returned an error message during computation. This happens when the length of a segments gets negative, that is, $t_i < 0$ for some index $i$.
\subsection{Benchmark 2}
\label{subsec:B2}
Assume a non-linear system of the form~\cite[p. ~334.]{Khalil:2002}
\begin{align*}
\dot{x}_1(t) & = - x_2(t) + x_1(t)x_3(t)\,, \\
\dot{x}_2(t) & = x_1(t) + x_2(t)x_3(t)\,, \\
\dot{x}_3(t) & = -x_3(t) - x_1(t)^2 - x_2(t)^2 + x_3(t)^2\,.
\end{align*}
We will investigate the behaviour of our method in dependence on the number of segments $N$ for problem formulations \eqref{eq:len_all}, \eqref{eq:pen_match_reg} and \eqref{eq:gap_pen_reg_diff}. Similarly to the previous Benchmark \ref{subsec:B1} we put $c_I = [1, 1, 1] \in \RR^3$. All the results are in Tab. \ref{tab:B2BFGS}.
\begin{table}
\centering
\begin{tabular}{ll}
\begin{tabular}{l|*{6}{c}}
N & $5$  &  $10$  &  $15$  &  $20$  &  $25$  &  $30$ \\
NIT & $35$  &  $33$  &  $31$  &  $48$  &  $47$  &  $49$ \\
S & $1$  &  $1$  &  $1$  &  $1$  &   $3$  &  $3$ \\ 
\end{tabular}
&
\begin{tabular}{l|*{6}{c}}
 N & $5$  &  $10$  &  $15$  &  $20$  &  $25$  &  $30$ \\
 NIT & $18$  &  $23$  &  $28$  &  $37$  &  $41$ &  $51$ \\
 S & $1$  &  $1$  &  $3$  &  $3$  &  $3$  &  $3$ \\ 
\end{tabular} \\
 \\
\begin{tabular}{l|*{6}{c}}
 N & $5$  &  $10$  &  $15$  &  $20$  &  $25$  &  $30$ \\
 NIT & $83$  &  $37$  &  $55$  &  $56$  &  $56$ &  $76$ \\
 S & $1$  &  $1$  &  $1$  &  $1$  &  $1$  &  $3$ \\ 
\end{tabular}
&
\begin{tabular}{l|*{6}{c}}
 N & $5$  &  $10$  &  $15$  &  $20$  &  $25$  &  $30$ \\
 NIT & $18$  &  $31$  &  $45$  &  $53$  &  $66$ &  $65$ \\
 S & $1$  &  $1$  &  $1$  &  $1$  &  $1$  &  $1$ \\ 
\end{tabular} \\
\\
\begin{tabular}{l|*{6}{c}}
 N & $5$  &  $10$  &  $15$  &  $20$  &  $25$  &  $30$ \\
 NIT & $23$  &  $229$  &  $269$  &  $400$  &  $400$  &  $400$ \\
 S & $1$  &  F  &  F  &  F  &  F  &  F \\
\end{tabular}
&
\begin{tabular}{l|*{6}{c}}
 N & $5$  &  $10$  &  $15$  &  $20$  &  $25$  &  $30$ \\
 NIT & $33$  &  $64$  &  $76$  &  $126$  &  $171$  &  $228$ \\
 S & F  &  F  &  F  &  F  &  F  &  F \\     
\end{tabular}
\end{tabular}
\caption{Benchmark 2: Block diagonal and banded BFGS approximations in the left and the full BFGS approximation in the right. There are results for formulation \eqref{eq:len_all} in the top; for  \eqref{eq:pen_match_reg} in the middle and for \eqref{eq:gap_pen_reg_diff} in the bottom. }
\label{tab:B2BFGS}
\end{table}

One can observe that problem formulation \eqref{eq:len_all} requires the least number of iterations. However, in some cases the method terminated because of the minimum step-length was reached. Problem formulation \eqref{eq:pen_match_reg} needs more iterations to finish, however, one again we obtained the desired solution with no fail attempts. Contrary to this the problem formulation \eqref{eq:gap_pen_reg_diff} yields the poorest results. Whenever we tried to verify results by simulation the computed solution did not meet our criteria $\| x_0^1 - c_I \|_{E_I} \leq 1 + 10^{-4}$ and $\| \Phi(\sum_{i = 1}^N t_i,x_0^1) - c_U \|_{E_U} \leq 1 + 10^{-4}$.
\subsection{Benchmark 3}
\label{subsec:B3}
In the end let us compare these three different problem formulations on a linear system. Assume we have the dynamics given by
\[
\dot{x}(t) = Ax(t)\,,
\]
where matrix $A$ is the same as in Benchmark 1 in \ref{subsec:B1}. We set $c_I = [1, \ldots, 1]^T \in \RR^n$. The results are in tables \ref{tab:B3PF03BFGS}, \ref{tab:B3PF04BFGS} and \ref{tab:B3PF09BFGS}.
\begin{table}
\centering
\begin{tabular}{*{24}{c}}
 n & N & NIT & S & n & N & NIT & S & n & N & NIT & S & n & N & NIT & S \\
\hline
 $10$  &  $5$  &  $28$  &  $1$  &  $20$  &  $5$  &  $33$  &  $1$  &  $30$  &  $5$  &  $30$ &  $1$  &  $40$  &  $5$  &  $34$  &  $1$ \\
 &  $10$  &  $31$  &  $1$  &    & $10$  &  $39$  &  $1$  &    &  $10$  &  $172$  &  $1$  &   &  $10$  &  $29$ &  $1$ \\
 &  $15$  &  $400$  &  $2$  &    &  $15$  &  $36$  &  $1$  &   &  $15$  &  $39$  &  $1$  &   &  $15$  &  $37$  &  $1$ \\
 &  $20$  &  $48$ &  $1$  &  &  $20$  &  $41$  &  $1$  &   &  $20$  &  $172$  &  $1$  &   &  $20$  &  $32$  &  $1$ \\
 &  $25$  &  $42$  &  $1$  &  &  $25$  &  $36$  &  $1$  &   &  $25$  &  $108$  &  $1$  &   &  $25$  &  $87$  &  $3$ \\
 &  $30$  &  $35$  &  $1$  &  &  $30$  &  $39$  &  $1$  &   &  $30$  &  $44$  &  $1$  &    &  $30$  &  $51$  &  $1$ \\  [5mm]
n & N & NIT & S & n & N & NIT & S & n & N & NIT & S & n & N & NIT & S \\
\hline
 $10$  &  $5$  &  $28$  &  $1$  &  $20$  &  $5$  &  $36$  &  $1$  &  $30$  &  $5$  &  $33$ &  $1$  &  $40$  &  $5$  &  $34$  &  $1$ \\
  &  $10$  &  $38$  &  $1$  &    & $10$  &  $46$  &  $1$  &  &  $10$  &  $54$  &  $1$  &   &  $10$  &  $56$ &  $1$ \\
  &  $15$  &  $44$  &  $3$  &    &  $15$  &  $54$  &  $1$  &  & $15$  &  $70$  &  $1$  &   &  $15$  &  $66$  &  $1$ \\
  &  $20$  &  $81$  & $1$  &   &  $20$  &  $81$  &  $1$  &  &  $20$  &  $50$  &  $3$  &  &  $20$  &  $50$  &  $3$ \\
  &  $25$  &  $69$  &  $3$  &    &  $25$  &  $121$ &  $1$  &  &  $25$  &  $134$  &  $1$  &    &  $25$  &  $159$  &  $1$ \\
 &  $30$  &  $78$  &  $1$  &   &  $30$  &  $100$  &  $1$  &  &  $30$  &  $116$  &  $1$  &    &  $30$  &  $133$  &  $1$ \\  
\end{tabular}
\caption{Benchmark 3 problem formulation \eqref{eq:len_all}: there are results for block diagonal BFGS approximation in the top and the full BFGS approximation in the bottom. }
\label{tab:B3PF03BFGS}
\end{table}
\begin{table}
\centering
\begin{tabular}{*{24}{c}}
 n & N & NIT & S & n & N & NIT & S & n & N & NIT & S & n & N & NIT & S \\
\hline
 $10$  &  $5$  &  $28$  &  $1$  &  $20$  &  $5$  &  $400$  &  F        &  $30$  &  $5$  &  $43$ &  $1$  &  $40$  &  $5$  &  $124$  &  F \\
  &  $10$  &  $30$  &  $1$  & &  $10$  &  $42$  &  $1$  &   &  $10$  &  $200$  &  F     &    &  $10$  &  $59$  & $1$ \\
  &  $15$  &  $50$  &  $1$  & &  $15$  &  $55$  &  $1$  &   &  $15$  &  $54$  &  $1$  &   &  $15$  &  $54$  &  $1$ \\
  &  $20$  &  $45$  &  $1$ & &  $20$  &  $49$  &  $1$  &   &  $20$  &  $54$  &  $1$  &   &  $20$ &  $53$  &  $1$ \\
  &  $25$  &  $44$  &  $1$  & &  $25$  &  $54$  &  $1$  &  &  $25$  &  $59$  &  $1$  &   &  $25$  &  $58$  &  $1$ \\
  &  $30$  &  $55$  &  $1$  & &  $30$  &  $52$  &  $1$  &   &  $30$  &  $53$  &  $1$  &   &  $30$  &  $66$  &  $1$ \\[5mm]
n & N & NIT & S & n & N & NIT & S & n & N & NIT & S & n & N & NIT & S \\
\hline
 $10$  &  $5$  &  $33$  &  $1$  &  $20$  &  $5$  &  $110$  &  F  &  $30$  &  $5$  &  $49$ &  $1$  &  $40$  &  $5$  &  $136$  &  F \\
   &  $10$  &  $55$  &  $1$  &   &  $10$  &  $63$  &  $1$  &  &  $10$  &  $232$  &  F  &  &  $10$  &  $400$  & $2$ \\
   &  $15$  &  $134$  &  $1$  &   &  $15$  &  $83$  &  $1$  &  &  $15$  &  $100$  &  $1$  &  &  $15$  &  $109$  &  $1$ \\
   &  $20$  &  $81$  &  $1$  &   &  $20$  &  $134$  &  $1$  &  &  $20$  &  $118$  &  $1$  &  & $20$  &  $400$  &  F \\
   &  $25$  &  $101$  &  $1$  &   &  $25$  &  $149$  &  $1$  &  &  $25$  &  $141$  &  $1$  &  &  $25$  &  $182$  &  $1$ \\
   &  $30$  &  $100$  &  $1$  &   &  $30$  &  $132$  &  $1$  &  &  $30$  &  $167$  & $1$  &  &  $30$  &  $177$  &  $1$ \\
\end{tabular}
\caption{Benchmark 3 problem formulation \eqref{eq:pen_match_reg}: there are results for block diagonal BFGS approximation in the top and the full BFGS approximation in the bottom. }
\label{tab:B3PF04BFGS}
\end{table}
\begin{table}
\centering
\begin{tabular}{*{24}{c}}
 n & N & NIT & S & n & N & NIT & S & n & N & NIT & S & n & N & NIT & S \\
\hline
 $10$  &  $5$  &  $393$  &  F       &  $20$  &  $5$  &  $400$  &  F        &  $30$  &  $5$  &  $400$ &  F       &  $40$  &  $5$  &  $400$  &  F \\
   &  $10$  &  $400$  &  F       &   &  $10$  &  $400$  &  F       &    &  $10$  &  $400$  &  F       &    &  $10$  &  $400$  &  F \\
   &  $15$  &  $400$  &  $2$  &    &  $15$  &  $400$  &  F       &   &  $15$  &  $400$  &  F        &    &  $15$  &  $400$  &  F \\
   &  $20$  &  $400$  &  $2$  &    &  $20$  &  $400$  &  F       &   &  $20$  &  $400$  &  F         &   &  $20$  &  $400$  &  F \\
   &  $25$  &  $400$  &  $2$  &    &  $25$  &  $400$  &  $2$  &    &  $25$  &  $400$  &  $2$  &    &  $25$  &  $400$  &  F \\
  &  $30$  &  $400$  &  F       &    &  $30$  &  $400$  &  F      &    &  $30$  &  $400$  &  F          &    &  $30$  &  $400$  &  F \\[5mm]
n & N & NIT & S & n & N & NIT & S & n & N & NIT & S & n & N & NIT & S \\
\hline
 $10$  &  $5$  &  $52$  &  F      &  $20$  &  $5$  &  $118$  &  $1$  &  $30$  &  $5$  &  $137$ &  F             &  $40$  &  $5$  &  $198$  &  F \\
   &  $10$  &  $296$  &  F     &    & $10$  &  $400$  &  F        &    &  $10$  &  $400$  &  F       &    &  $10$  &  $400$  &  F \\
   &  $15$  &  $163$  &  F      &    &  $15$  &  $400$  &  F       &   &  $15$  &  $400$  &  F        &    &  $15$  &  $400$  &  F \\
   &  $20$  &  $209$ &  F     &    &  $20$  &  $400$  &  F       &    &  $20$  &  $400$  &  F         &    &  $20$  &  $400$  &  F \\
   &  $25$  &  $292$  &  F     &    &  $25$  &  $364$  &  $1$  &   &  $25$  &  $362$  &  $1$  &    &  $25$  &  $400$  &  F \\
  &  $30$  &  $314$  &  $1$  &   &  $30$  &  $392$  &  $1$  &   &  $30$  &  $400$  &  F     &    &  $30$  &  $400$  &  F \\
\end{tabular}
\caption{Benchmark 3 problem formulation \eqref{eq:gap_pen_reg_diff} there are results for block diagonal BFGS approximation in the top and the full BFGS approximation in the bottom. }
\label{tab:B3PF09BFGS}
\end{table}

We can see in Table \ref{tab:B3PF03BFGS} that our method found an error trajectory in all setups for problem formulation \eqref{eq:len_all}. Problem formulation \eqref{eq:pen_match_reg}, with results in Tab. \ref{tab:B3PF04BFGS}, performed well with only a few failed attempts. In the end, problem formulation \eqref{eq:gap_pen_reg_diff} failed many times as you can see in Tab. \ref{tab:B3PF09BFGS}. To this end problem formulation \eqref{eq:len_all} can be said to be superior to \eqref{eq:pen_match_reg} and \eqref{eq:gap_pen_reg_diff} since it produces better results on all three benchmarks.
\subsection{Trust-region SQP}
\label{subsec:TRSQP}
An alternative approach to line search SQP is trust-region SQP \cite[Alg.~18.4]{Nocedal:2006}. In order to check the performance of the trust-region method we recomputed all benchmarks for problem formulation \eqref{eq:len_all} with the block-diagonal BFGS approximation and received similar results for Benchmarks \ref{subsec:B1} and \ref{subsec:B2}, and worse results measured in terms of iterations for Benchmark \ref{subsec:B3}. When one uses BFGS for the approximation of the Hessian then line search SQP is performing well enough and trust-region SQP does not bring any improvement. The reason behind choosing problem formulation \eqref{eq:len_all} for the comparison is that it performs best for the experiments shown above. 

We also tried trust-region SQP for the other problem formulations \eqref{eq:pen_match_reg} and \eqref{eq:gap_pen_reg_diff}, and it does not show any improvement. Furthermore, for those problem formulations it is difficult to set a good maximal trust-region radius for trust-region SQP to converge.
\section{Conclusion}
\label{sec:Conclusion}
We presented a solution to the problem of finding an error trajectory of ordinary differential equations.  We considered several different constrained minimization problem formulations, that we solve using SQP. We discussed the influence of the structure of the formulation on the solution algorithm, and performed computational experiments that showed that Formulation~\eqref{eq:len_all} results in the the most successful method. Here, the KKT system features a block-diagonal Hessian matrix and sparse Jacobian of constraints. 

As such, this paper gives insight into optimization techniques for solving the falsification problem of dynamical systems where, up to now, optimization techniques have only been applied as a black box~\cite{Abbas:2011,S-TaLiRo:2011,KuratkoRatschan:2014,Nghiem:2010,Zutshi:2013}. 
\section{Acknowledgement}
\label{sec:Thanks}
We would like to thank Ladislav Luk\v{s}an for his insightful tips about the implementation of various updating schemes for the Hessian matrix. Moreover, we would like to thank Miroslav Rozlo\v{z}n\'{i}k with whom we discussed the structure of the saddle point matrix and the solution method we can apply.
\bibliographystyle{abbrv}
\bibliography{../bibliography/kuratko}

\section*{Appendix}
\label{sec:Appendix}
First we address the linear independence constraint qualification (LICQ) condition for our vectors of constraints.
\begin{lemma}
\label{lem:grad_pen_match}
Let $A \in \RR^{N(n+1)\times (N-1)n}$ be a matrix of the form
\[
A = \begin{bmatrix}
M_1 &  & & & \\
v_1^T &   & & & \\
I & M_2 & & & \\
 & v_2^T & & & \\
  &  I & \ddots &  & \\
 &   & \ddots & v_{N-2}^T & \\
 &  & & I & M_{N-1} \\
 &  & &  & v_{N-1}^T \\
 &  & &  & I \\ 
  &  & &  & 0 \\ 
\end{bmatrix}\ ,
\]
where for $1 \leq i \leq N-1$, $M_i \in \RR^{n \times n}$, $v_i \in \RR^n$, $0 \in \RR^{1 \times n}$, and $I \in \RR^{n \times n}$ is the identity matrix. Then the matrix $A$ has full-column rank.
\begin{proof}
We prove this Lemma by contradiction. Suppose matrix $A$ does not have full-full column rank. Then there exists a non-zero vector $x \in \RR^{(N-1)n}$ such that $x = \left[ x_1, \ldots, x_{N-1}\right]^T$, where $x_i \in \RR^n$ for $1 \leq i \leq N-1$, for which $Ax = 0$. Rewriting the equation $Ax = 0$, we get
\begin{align*}
x_i + M_{i+1}x_{i+1} & = 0\,,   \\
x_{N-1} & = 0\,, 
\end{align*}
for $1 \leq i \leq N-2$. We can observe that substituting backwards from $x_{N-1} = 0$ we get all $x_i = 0$, $1 \leq i \leq N-1$.  This is a contradiction with our assumption that $x \in \RR^{(N-1)n}$ is a non-zero vector.
\end{proof}
\end{lemma}
\begin{lemma}
\label{lem:grad_len_all}
Let $A \in \RR^{N(n+1)\times ((N-1)n+2)}$ be a matrix of the form
\[
A = \begin{bmatrix}
w_1 & M_1 & &   & & & \\
&v_1^T &   &&  & & \\
&I &  M_2 & & & & \\
 & & v_2^T & & & & \\
 & &  I & \ddots & & & \\
& &   & \ddots & v_{N-2}^T &  & \\
& &  & & I & M_{N-1}  & \\
& &  & &  & v_{N-1}^T & \\
& &  & &  & I & w_2\\ 
 & &  & &  & 0 & \omega\\ 
\end{bmatrix}\ ,
\]
where for $1 \leq i \leq N-1$, $M_i \in \RR^{n \times n}$, vectors $v_i, w_1, w_2 \in \RR^n$, $\omega \in \RR$, and $I \in \RR^{n \times n}$ is the identity matrix. If $w_1 \in \RR^n$ is a non-zero vector and $\omega \neq 0$, then the matrix $A$ has full-column rank.
\begin{proof}
We prove this Lemma by contradiction. Suppose columns in $A$ are linearly dependent, therefore, there exists a non-zero vector $x = \left[ \alpha, x_1, \ldots, x_{N-1}, \beta \right]^T \in \RR^{(N-1)n+2}$ with $x_i \in \RR^n$, $1 \leq i \leq N-1$ $\alpha \in \RR$, and $\beta \in \RR$ so that
\begin{align*}
\alpha w_1 + M_1x_1 & = 0\,,\\
x_i + M_{i+1}x_{i+1} & = 0\,,   \\
x_{N-1} + \beta w_2 & = 0\,, \\ 
\beta\omega = 0 \,,
\end{align*}
for $1 \leq i \leq N-2$. Since we assume $\omega$ to be a non-zero scalar, therefore, we get $\beta = 0$. It follows that $x_{N-1} = 0 \in \RR^n$. If we substitute into formulae above we obtain $x_i = 0 \in \RR^n$ for $1 \leq i \leq N-2$. Therefore, also $\alpha w_1 = 0 \in \RR^{n}$. For $\alpha \neq 0$ this is only possible if $w_1 = 0 \in \RR^{n}$. This is contradiction with the assumption that $w_1$ is a non-zero vector. 
\end{proof}
\end{lemma}
\end{document}